\documentclass{article}
\usepackage{graphicx} 
\usepackage{amssymb}
\usepackage{amsmath, geometry}
\usepackage{hyperref}
\title{Discrete Non-Abelian X-Ray Transforms}
\author{Pranav Gupta and Roman G. Novikov}
\date{}
\begin{document}

\maketitle
\begin{abstract}
    We define a discrete version of the non-abelian X-ray transform, going back in particular to Manakov, Zakharov (1981) and Strichartz (1982). We extend to this transform non-overdetermined reconstruction results obtained for the abelian case in the recent article by Novikov, Sharma  (2025). In addition, we establish relations with the continuous non-abelian X-ray transform. In this respect, our results include an explicit and exact non-overdetermined layer-stripping reconstruction procedure for piecewise constant matrix-valued functions from their continuous non-abelian X-ray transform. To our knowledge, this result is new even for the classical X-ray transform.
\end{abstract}
\begin{center}
\begin{minipage}{0.9\textwidth} 
\end{minipage}
\end{center}

\section{Introduction}
Let $f$ be a function on $\mathbb{R}^{d}$ that takes values in $M(n,\mathbb{C})$, that is in $n \times n$ complex matrices, and let $T$ denote the set of all oriented straight lines $\gamma$ in $\mathbb{R}^d$. \\
We consider the non-abelian X-ray transform $\mathcal{S}$ of $f$:
\begin{equation}
\mathcal{S}(f)(\gamma) = \lim_{s \rightarrow+\infty}\psi^{+}_{\gamma}(s), \ \gamma \in T, \ s \in \mathbb{R},
\label{1.1}
\end{equation}
where $\psi^{+}_{\gamma}(s)$ is the solution to: \\
\begin{equation}
    \frac{d}{ds} \psi_{\gamma}^{+}(s) = f(x_{\gamma}(s))\psi^{+}_{\gamma}(s),
\label{1.2}
\end{equation}
\begin{center}
    $\lim_{s\rightarrow-\infty}\psi^{+}_{\gamma}(s) = \mathbf{I},$
\end{center}
where $x_{\gamma}(s)$ denotes a parameterization of $\gamma$ by natural parameter $s,$ $\mathbf{I}$ is the identity matrix. \\
Here, we assume that $f$ is sufficiently regular and decays sufficiently fast as $|x| \rightarrow \infty$.\\
If $f$ takes values in $\mathbb{C},$ then
$\mathcal{S}(f) = \exp(\mathcal{P}f)$, where $\mathcal{P}f$ is the classical X-ray transform defined by
\begin{equation}
\mathcal{P}f(\gamma) = \int_{\gamma}f(y)dy,
\label{1.3}
\end{equation}
and arising, in particular, in the X-ray tomography; see, for example, \cite{N2001}. \\\\
In the matrix case, in connection with the definition of $\mathcal{S}(f)$, methods for finding $f$ from $\mathcal{S}(f)$, generalizations and applications, see \cite{MZ1981}, \cite{V1992}, \cite{N2002}, \cite{E2004}, \cite{N2019}, \cite{MNP2021}, \cite{PSU2023}, and references therein.
\\\\
In the present work, we define and study a discrete version of the non-abelian X-ray transform $\mathcal{S}(f)$, which reduces in the abelian case to the discrete X-ray transform going back to \cite{S1982}). \\
Let $F$ be a function on $\mathbb{Z}^d$ that takes values in $GL(n,\mathbb{C})$, that is in the group of $n \times n$ complex matrices with non-zero determinant. We consider the discrete non-abelian X-ray transform $\mathcal{S}$ of $F$:
\begin{equation}
    \mathcal{S}(F)(\gamma) = \prod^{\longleftarrow}_{y \in \gamma \cap \mathbb{Z}^d}F(y), \ \gamma \in T',
    \label{1.4}
\end{equation}
\begin{equation}
    T' = \bigcup_{\zeta \in \mathbb{Z}^{d}}T_{\zeta}, \ T_{\zeta} := \{\gamma \in T: \zeta \in \gamma\}.
\label{1.5}
\end{equation}
Note that in \eqref{1.4} , the matrix product is considered in the following sense
\begin{equation}
    \prod^{\longleftarrow}_{y \in \gamma \cap \mathbb{Z}^d}F(y) = ...F(y_j)F(y_{j-1})...F(y_{i+1})F(y_i)...
\label{1.6}
\end{equation}
where $\{y_j\}$ are the points of $\gamma \cap \mathbb{Z}^d$, numbered from right to left according to the orientation of $\gamma$, and we assume that $F(y) \rightarrow \mathbf{I}$ sufficiently fast as $|y| \rightarrow \infty$.
\\
If $F = \exp(f)$, where $f$ is a complex-valued function on $\mathbb{Z}^d$, then:
\begin{equation}
 \mathcal{S}(F)(\gamma) = \exp(\mathcal{P}f(\gamma)),
 \label{1.7}
\end{equation}
\begin{equation}
 \mathcal{P}f(\gamma) = \sum_{y \in\gamma \cap \mathbb{Z}^d}f(y), \ \gamma \in T'.
 \label{1.8}
\end{equation}
Here, $\mathcal{P}f$ is the discrete X-ray transform as in \cite{NS2025}. The definition of this transform goes back to \cite{S1982}, at least in dimension $d = 2$ and for straight lines $\gamma$ with rational directions (as in formula \eqref{1.15} below).\\
Our definition of the non-abelian discrete X-ray transform $\mathcal{S}$ in \eqref{1.4} is also motivated by the formula:
\begin{equation}
    \mathcal{S}(f)(\gamma) = \lim_{\Delta \rightarrow 0}...\exp(f(y_j)\Delta)\exp(f(y_{j-1})\Delta)...\exp(f(y_{i+1})\Delta)\exp(f(y_{i})\Delta)...
\label{1.9}
\end{equation}
\begin{center}
    $\gamma = \bigcup_{i =-\infty}^{i = \infty} \Lambda_{i}$, \ $\Lambda_{i} = [y_i, y_{i+1})$, $|y_{i+1} - y_i| = \Delta$,
\end{center}
where $\mathcal{S}(f)$ is defined as in \eqref{1.1}, \eqref{1.2}, and $f$ is continuous on $\mathbb{R}^d$ with sufficient decay at infinity. 
For more information on the multiplicative integral arising in \eqref{1.9}, see for example, \cite{MZ1981}, and \cite{S2007}.\\\\
The transform $\mathcal{S}$ defined in \eqref{1.4} is not the only natural version of discrete non-abelian X-ray transforms. A different natural version $\mathcal{S}^{\star}$ of such transforms is defined in formula \eqref{3.2.8} in Section 5. One more natural version of such transforms, especially in connection with the theory of solitons, can be defined proceeding from \cite{MZ1981} and \cite{FN1991}.
\\\\
In the present work, we extend results in \cite{NS2025} on finding $f$ from $\mathcal{P}f$ defined by \eqref{1.8} to the non-abelian case, that is to the problem of finding $F$ from $\mathcal{S}(F)$ defined by \eqref{1.4}. \\
Our basic assumptions on $F$ are as follows:
\begin{equation}
    \det(F) \neq 0,
\label{1.10}
\end{equation}
\begin{equation}
supp \ (F - \mathbf{I{)}} \subseteq B_r \cap \mathbb{Z}^d,
\label{1.11}
\end{equation}
where
\begin{equation}
    B_r := \{x \in \mathbb{R}^d: |x| \leq r\}, \ r > 0.
\label{1.12}
\end{equation}
Let
\begin{equation}
    N_r := \#(B_r \cap \mathbb{Z}^d ),
\label{1.13}
\end{equation}
where $\#$ denotes cardinality. \\
Let
\begin{equation}
    \Omega := \{\theta \in \mathbb{S}^{d-1}: \theta = z/|z|, \ z \in \mathbb{Z}^d \setminus\{0\}\}, 
\label{1.14}
\end{equation}
\begin{equation}
    \mathcal{T} := \{\gamma \in T': \hat{\gamma} \in \Omega \},
\label{1.15}
\end{equation}
where $\hat{\gamma}$ denotes the direction of $\gamma$, and the sets $\Omega$ and $\mathcal{T}$ are considered as the sets of rational directions and rational rays, respectively.\\
The results of the present work can be summarized as follows:
\begin{itemize}
    \item Under assumptions \eqref{1.10}, \eqref{1.11}, we give non-overdetermined reconstructions of $F$ from $\mathcal{S}(F)$ restricted to an appropriate $\Gamma_r \subseteq T'$ with $\#\Gamma_r = N_r$. The simplest case admitting irrational rays is given in Subsection 2.1, whereas a layer-stripping reconstruction involving only rational rays is given in Subsection 2.2.
    \item We present the relations between the scalar weighted discrete X-ray transfrom $\mathcal{P}_W$, considered in \cite{NS2025}, and the discrete non-abelian X-ray transform $\mathcal{S}$; see Section 3. In particular, Theorem 3.1 shows that the transform $\mathcal{P}_W(u)$ arises as reduction of  $\mathcal{S}(F_\theta)$, where $F_\theta$ is a matrix-valued function. This result serves as a discrete analogue of well-known results in the continuous case; see \cite{N2002}, \cite{N2019}.
    \item We analyze relations between continuous and discrete non-abelian X-ray transforms: \\
    In Section~4, we establish relations between the discrete non abelian X-ray transform defined in \eqref{1.4} and the continuous non-abelian X-ray transform defined by \eqref{1.1}, \eqref{1.2} for the case when $f$ is a sum of regularized delta functions, where the initial delta functions are supported on $B_r \cap \mathbb{Z}^d.$ \\
    In Section 5, we suggest an additional natural version $\mathcal{S}^{\star}$ of discrete non-abelian X-ray transform arising from the continuous non-abelian X-ray transform $\mathcal{S}(f)$ of piecewise constant matrix function $f$; see formula \eqref{3.2.8}. For the transform $\mathcal{S}^{\star}$, we give a non-overdetermined layer-stripping reconstruction in Theorem 5.5, thus yielding an exact non-overdetermined layer-stripping reconstruction of a piecewise constant function $f$ from its continuous non-abelian X-ray transform $\mathcal{S}(f)$. To our knowledge, this reconstruction is new, even for the case of classical X-ray tomography. 
\end{itemize}
Numerical implementation of the aforementioned reconstructions will be studied elsewhere.

\section{Non-Overdetermined Reconstructions}
In this section we give formulas for non-overdetermined reconstructions of $F$ from its discrete non-abelian X-ray transform $\mathcal{S}(F)$. These formulas are similar to those in \cite{NS2025} for the abelian case. Nevertheless, an essential non-abelian structure arises in formula \eqref{2.30} used in the layer-stripping reconstruction of $F$ from $\mathcal{S}(F)$ given on rational rays. Finally, in this section, we formulate two open questions based on related conjectures for the abelian case. 
\subsection{Irrational Rays}
The simplest reconstruction of $F$ on $\mathbb{Z}^d$ from $\mathcal{S}(F)$ on $T'$ defined by \eqref{1.5} consists of the following formula:
\begin{equation}
    F(x) = \mathcal{S}(F)(\gamma_{x,\theta(x)}), \ x \in \mathbb{Z}^d, \ \theta(x) \in \mathbb{S}^{d-1} \setminus \Omega,
\label{2.16}
\end{equation}
where $\gamma_{x,\theta(x)}$ is the straight line passing through $x$ with direction $\theta$ = $\theta(x)$, and $\Omega$ is defined by \eqref{1.14}. One can see that:
\begin{itemize}
    \item 
    Formula \eqref{2.16} involves $\mathcal{S}(F)$ only on $\Gamma$, where $\Gamma$ is defined by:
    \begin{equation}
    \Gamma := \{\gamma = \gamma_{x,\theta(x)}: x \in \mathbb{Z}^d, \ \theta(x) \in \mathbb{S}^{d-1} \setminus \Omega \}.
    \label{2.17}
    \end{equation}
    \item
    The elements of $\Gamma$ are in one-to-one correspondence with those of $\mathbb{Z}^d$.
    \item
    Formula \eqref{2.16} holds for arbitrary $GL(n,\mathbb{C})-$valued function on $\mathbb{Z}^d$ (without assumption \eqref{1.11}).
    \item 
    Under the additional assumption \eqref{1.11}, formula \eqref{2.16} reconstructs $F$ on $B_r \cap \mathbb{Z}^d$ from $\mathcal{S}(F)$ on $\Gamma_r$ only, where:
    \begin{equation}
    \Gamma_r := \{\gamma = \gamma_{x,\theta(x)}: x \in B_r \cap \mathbb{Z}^d, \ \theta(x) \in \mathbb{S}^{d-1} \setminus \Omega \}, \ \ \#\Gamma_r = N_r,
    \label{2.18}
    \end{equation}
    where $N_r$ is defined by \eqref{1.13},
    \item
    Under assumptions \eqref{1.10}, \eqref{1.11}, formula \eqref{2.16} also holds when $\theta(x) \in \mathbb{S}^{d-1} \setminus \Omega_{\rho}, \ \rho > 2r$, where 
    \begin{equation}
    \Omega_{\rho} := \{\theta: \theta = z/|z|, z = (z_1,...,z_d) \in \mathbb{Z}^d, |z|/\gcd\{z_1,...,z_d\} \leq \rho \}, \ \rho \geq 1.
    \label{2.19}
    \end{equation}
    In addition, in formula \eqref{2.18}, $\Omega$ can be replaced by $\Omega_{\rho}$, where $\rho > 2r$.
\end{itemize}
Thus, formula \eqref{2.16} gives a non-overdetermined reconstruction of $F$ from $\mathcal{S}(F)$ for the case of irrational rays, and also for the case of rays with sufficient degree of irrationality under assumption \eqref{1.11}.

\subsection{Rational Rays}
We consider the set $\mathcal{T}^{\star} \subset T' $ of rational rays defined by the following formulas:
\begin{equation}
    \mathcal{T}^{\star} := \bigcup_{z \in \mathbb{Z}^d} \gamma_z,
\label{2.20}
\end{equation}
where
\begin{equation}
\gamma_z = \{y \in \mathbb{R}^d: y = z + t\hat{z}^{\star}, t\in \mathbb{R}\}, \ z \in \mathbb{Z}^d,\  z_1^2 + z_2^2 \neq 0,
\label{2.21}
\end{equation}
\begin{equation}
    \hat{z}^{\star} := \frac{-z_2e_1 + z_1e_2}{\sqrt{z_1^2 + z_2^2}}, \ z = \sum_{i = 1}^d z_ie_i,
\label{2.22}
\end{equation}
and
\begin{equation}
    \gamma_z = \{y \in \mathbb{R}^d: y = z + t\theta, t\in \mathbb{R}\}, \ \theta \in \Omega, \ \theta = (\theta_1,\theta_2,0,...,0), \ z_1^2 + z_2^2 = 0, 
\label{2.23}
\end{equation}
where $e_1,...,e_d$ are the standard basis vectors in $\mathbb{R}^d$, $\Omega$ is defined by \eqref{1.14}. We also consider:

\begin{equation}
    \mathcal{A}_{\alpha,\beta} := \{z \in \mathbb{Z}^d: \alpha \leq \sqrt{z_1^2 + z_2^2} \leq \beta\},
\label{2.24}
\end{equation}
\begin{equation}
    \mathcal{T}^{\star}_{\alpha,\beta} := \{\gamma_z \in \mathcal{T}^{\star}: \alpha \leq \sqrt{z_1^2 + z_2^2} \leq \beta\}, \  0 \leq \alpha \leq \beta.
\label{2.25}
\end{equation}
\\
\textbf{Theorem 2.1} \textit{Let $F$ satisfy \eqref{1.10}, \eqref{1.11}, then $\mathcal{S}(F)$ on $\mathcal{T}^{\star}$ uniquely determines $F$, where $\mathcal{T}^{\star}$ is defined by \eqref{2.20}. Moreover, this determination is given via the two-dimensional formulas \eqref{2.29}, \eqref{2.30} for $d = 2$, which are adopted for each slice with fixed $z_3,...,z_d$, for $d \geq 3$. In addition, in this framework, $\mathcal{S}(F)$ on $\mathcal{T}^{\star}_{\alpha, \beta}$, where $\beta \geq r$, uniquely determines $F$ on $\mathcal{A}_{\alpha,\beta}$.}
\\\\
\textbf{Proof of Theorem 2.1} First, we consider the case $d = 2$. Let:
\begin{equation}
    S_1 := \{z \in B_r \cap \mathbb{Z}^2: |z| \geq |{\zeta}|, \ \forall \zeta \in B_r \cap \mathbb{Z}^2 \}, \ S_J = \{0\},
\label{2.26}
\end{equation}
\begin{equation}
    S_{i+1} := \{z \in B_r \cap \mathbb{Z}^2\setminus \cup_{k = 1}^i S_k: |z| \geq |{\zeta}|, \ \forall \zeta \in (B_r \cap \mathbb{Z}^2) \setminus (\cup_{k = 1}^i S_k) \}, \ i+1 \leq J,
\label{2.27}
\end{equation}
where $B_r$ is the disc defined in \eqref{1.12} for $d = 2$. Recall that 
\begin{equation}
    B_r \cap \mathbb{Z}^2 = \cup_{i = 1}^J S_i,\  \ S_i \cap S_j = \emptyset \textrm{ if } i \neq j.
\label{2.28}
\end{equation}
The reconstruction formulas are as follows:
\begin{equation}
    F(z) = \mathcal{S}(F)(\gamma_z), \ z \in S_1,
\label{2.29}
\end{equation}
\begin{equation}
    F(z) = (\prod^{\longleftarrow}_{\zeta \in {\gamma_z^+}\cap \atop(\cup_{k =1}^{i}S_k)}F(\zeta))^{-1}\mathcal{S}(F)(\gamma_z)(\prod^{\longleftarrow}_{\zeta \in {\gamma_z^-} \cap \atop(\cup_{k =1}^{i}S_k)}F(\zeta))^{-1}, \ z \in S_{i+1}, \ i = 1,2,...,J-1,
\label{2.30}
\end{equation} 
\begin{equation}
    {\gamma_z^+} := \{y \in \gamma_z: (y-z).\hat{\gamma}_z > 0 \},
    \atop
    {\gamma_z^-} := \{y \in \gamma_z: (y-z).\hat{\gamma}_z < 0 \},
\label{2.31}
\end{equation}
where the products are defined as the product in \eqref{2.16}. \\
More precisely, first, we reconstruct $F$ on $S_1$ by \eqref{2.29}. Then, assuming that $F$ is reconstructed on $\cup_{k =1}^i S_k$, we reconstruct $F$ on $S_{i+1}$ by \eqref{2.30}. In view of \eqref{2.28}, formulas \eqref{2.29}, \eqref{2.30} yield a reconstruction of $F$ on $B_r \cap \mathbb{Z}^2$ in $J$ inductive steps.
\\
In addition, formulas \eqref{2.29}, \eqref{2.30} yield a reconstruction of $F$ on $\mathcal{A}_{\alpha,\beta}$ from $\mathcal{S}(F)$ given on $\mathcal{T}^{\star}_{\alpha,\beta}$, for $\beta \geq r.$
\\
To get formulas \eqref{2.29}, \eqref{2.30}, we use that:
\begin{equation}
    \gamma_z = \gamma_z^{-} \cup \{z\} \cup \gamma_z^{+}, \ z \in \mathbb{Z}^2,
\label{2.32}
\end{equation}
\begin{equation}
    \gamma_z \cap (\cup_{k = i + 1}^{J} S_k) = \{z\}, \ z \in S_{i+1}, 
\label{2.33}
\end{equation}
and, that, under assumption \eqref{1.11}, we have
\begin{equation}
    \mathcal{S}(F)(\gamma_z) = (\prod^{\longleftarrow}_{\zeta \in {\gamma_z^+}\cap \atop\bigcup_{k =1}^{i}S_k}F(\zeta))F(z)(\prod^{\longleftarrow}_{\zeta \in {\gamma_z^-} \cap \atop\bigcup_{k =1}^{i}S_k}F(\zeta)), \ z \in S_{i+1}, 
\label{2.34}
\end{equation}
where $i = 0,1,...,J-1$, and the products in \eqref{2.34} are defined as $\mathbf{I}$ for $i = 0$. In view of \eqref{1.10}, formula \eqref{2.30} follows from \eqref{2.34}. \\
This completes the proof of Theorem 2.1 for $d = 2$.
\\
In dimension $d \geq 3$, we adopt the reconstruction formulas \eqref{2.29}, \eqref{2.30} for each two-dimensional plane 
\begin{equation}
    \Xi = \{x \in \mathbb{R}^d: x = (x_1,x_2,z_3,...,z_d), (x_1, x_2) \in \mathbb{R}^{d}\}, \  z_3,...,z_d \in \mathbb{Z},
\label{2.n.1}
\end{equation} where the origin $O$ of $\Xi$ is considered to be $x_0 = (0,0,z_3,...,z_d).$ This yields a proof of Theorem 2.1 for $d \geq 3$.
\\\\
\textbf{Remark 2.2}
\textit{If $F$ takes values in $\mathbb{C}\setminus\{0\}$, Theorem 2.1 reduces to the corresponding result in \cite{NS2025} for the case of abelian discrete X-ray transform, in view of formulas \eqref{1.7}, \eqref{1.8}.}
\\\\
\textbf{Remark 2.3}
\textit{In a similar way to the considerations in \cite{NS2025}, Theorem 2.1 can be extended to the case of rays parallel to the plane 
\begin{equation}
    \Xi_{a,b} := \textrm{span}\{a,b\},
\label{2.35}
\end{equation}
where $a,b \in \mathbb{Z}^d,$ and are linearly independent. In this case, the rays $\gamma_z$ are defined by the formulas:
\begin{equation}
    \gamma_z = \{y \in \mathbb{R}^d: y = z + tz^{\star}, \ t \in \mathbb{R}\}, \ z \in \mathbb{Z}^d, \ (z\cdot a)^2 + (z\cdot b)^2 \neq 0,
\label{2.36}
\end{equation}
\begin{equation}
    z^\star = \frac{-(z\cdot b)a + (z\cdot a)b}{|-(z\cdot a)b + (z\cdot b)a|},
\label{2.37}
\end{equation}
and
\begin{equation}
    \gamma_z = (\theta, z), \ \theta \in \Omega \cap \Xi_{a,b}, \ z \in \mathbb{Z}^d, \ (z \cdot a)^2 + (z \cdot b)^2 = 0.
\label{2.38}
\end{equation}
}
\\
Note that reconstructions considered in Theorem 2.1, Remarks 2.2 and 2.3 belong to the class of layer stripping reconstructions in inverse problems. In the domain of Radon-type transforms, such reconstructions go back to the Cormack reconstruction for the classical two dimensional Radon transform in \cite{C1963}. \\
In addition to Theorem 2.1, Remarks 2.2 and 2.3, there are the following open questions:
\\
\textbf{Question 1}
\textit{Whether $\mathcal{S}(F)$ on $\mathcal{T}^{\star}$ uniquely determines $F$ on $\mathbb{Z}^2$, where $F$ satisfies \eqref{1.10},}
\begin{equation}
    F(x) - \mathbf{I} = O(|x|^{-N}), \ x \in \mathbb{Z}^2, |x| \rightarrow \infty, \ \textrm{for some } N > 1,
\label{2.39}
\end{equation}
\textit{and $\mathcal{T}^{\star}$ is defined in \eqref{2.20} for $d = 2$?}
\\\\
\textbf{Question 2}
\textit{Whether $\mathcal{S}(F)$ on $\mathcal{T}^{\star}_{\alpha,\beta}$ uniquely determines $F$ on $\mathcal{A}_{\alpha,\beta}$, where $F$ satisfies \eqref{1.10},}
\begin{equation}
    F(x) - \mathbf{I} = O(|x|^{-\infty}), \ x \in \mathbb{Z}^2, |x| \rightarrow \infty,
\label{2.40}
\end{equation}
\textit{and $\mathcal{T}^{\star}_{\alpha,\beta}$ and $\mathcal{A}_{\alpha,\beta}$ are defined in \eqref{2.24},\eqref{2.25} for $d = 2$?} 
\\\\
Note that these questions for $F$ taking values in $\mathbb{C}\setminus\{0\}$ correspond to related conjectures in \cite{NS2025} for the abelian case.

\section{Reductions to Weighted Discrete X-Ray Transforms}
Recall that for continuous and discrete X-ray transforms \eqref{1.3} and \eqref{1.8}, their weighted versions are defined by the formulas \eqref{4.1} and \eqref{4.2}, respectively:
\begin{equation}
    \mathcal{P}_Wu(\gamma) = \int_{\gamma}W(y,\hat{\gamma})u(y)dy, \ \gamma \in T,
\label{4.1}
\end{equation}
\begin{equation}
    \mathcal{P}_Wu(\gamma) = \sum_{y \in \gamma \cap \mathbb{Z}^d}W(y,\hat{\gamma})u(y), \ \gamma \in T',
\label{4.2}
\end{equation}
where $W$ is a weight function, $T$ and $T'$ are as in formulas \eqref{1.1} and \eqref{1.5} , and $\hat{\gamma}$ is the direction of $\gamma.$ Usually, it is assumed that $W$ is real valued and strictly positive.\\
Recall that the transforms $\mathcal{P}_Wu$ in \eqref{4.1} can be seen as reductions of the continuous non-abelian X-ray transforms; see \cite{N2002} and \cite{N2019}. Discrete analogues of these results are given below as Theorem 3.1 and Remark 3.2:
\\\\
\textbf{Theorem 3.1}
\textit{Let 
\begin{equation}
    F_{\theta}(z) = \begin{pmatrix}
        w_{1,\theta}(z) & w_{2,\theta}(z)u(z) \\
        0 & 1
    \end{pmatrix}, \ z \in \mathbb{Z}^d, \ \theta \in \mathbb{S}^{d-1},
\label{4.1.1}
\end{equation}
where 
\begin{equation}
\begin{aligned}
    w_{1,\theta}, w_{2,\theta}, u \textrm{ are complex-valued, and  } w_{1,\theta}, w_{2,\theta} \textrm{ have no zeros,}
    \\
    \textrm{supp}(w_{1,\theta} -1), \ \textrm{supp}(w_{2,\theta}-1), \ \textrm{supp(u)} \subseteq B_r \cap \mathbb{Z}^d.
\end{aligned}
\label{4.1.2}
\end{equation}
Then, 
\begin{equation}
     \mathcal{S}(F_{\theta})(\gamma) = \begin{pmatrix}
        \mathcal{S}(w_{1,\theta})(\gamma) & P_Wu(\gamma) \\
        0 & 1
    \end{pmatrix}, \ \gamma \in T', \ \hat{\gamma} = \theta,
\label{4.1.3}
\end{equation}
where $\mathcal{S}(w_{1,\theta})(\gamma)$ is the discrete non-abelian X-ray transform of $w_{1,\theta}$, and $P_Wu(\gamma)$ is defined by \eqref{4.2}, where
\begin{equation}
    W(y_,\hat{\gamma}) = w_{2,\theta}(y)\prod_{z \in \gamma^+_{y,\theta}\cap\mathbb{Z}^d}w_{1,\theta}(z), \ \gamma \in T', \ \hat{\gamma} = \theta, \ y \in \gamma \cap \mathbb{Z}^d,
    \atop
    \gamma^{+}_{y, \theta} = \{x \in \gamma: (y - x)\cdot\theta > 0\}.
\label{4.1.4}
\end{equation}
}
\\\\
\textbf{Proof of Theorem 3.1}
Fix $\gamma \in T'$ such that $\hat{\gamma} = \theta,$ and let $\gamma \cap \mathbb{Z}^d \cap B_r = \{y_1,\cdots,y_n\},$ where $(y_i - y_j)\cdot\hat{\gamma} > 0$ if $i > j.$ Then, in view of \eqref{4.1.2},
\begin{equation}
    \mathcal{S}(F_{\theta})(\gamma) = F_\theta(y_n)\cdots F_\theta(y_1) \\
    = \begin{pmatrix}
        w_{1,\theta}(y_n) & w_{2,\theta}(y_n)u(y_n) \\
        0 & 1
    \end{pmatrix}
    \cdots
    \begin{pmatrix}
        w_{1,\theta}(y_1) & w_{2,\theta}(y_1)u(y_1) \\
        0 & 1
    \end{pmatrix},
\label{4.1.5}
\end{equation}
\begin{equation}
    \mathcal{P}_Wu(\gamma) = \sum_{k = 1}^n\bigg(w_{2,\theta}(y_k)\prod_{l = k+1}^{n}w_{1,\theta}(y_l)\bigg)u(y_k), \ \prod_{l = n+1}^{n}w_{1,\theta}(y_l) = 1,
\label{4.1.6}
\end{equation}
where $W$ is defined by \eqref{4.1.4}. It follows from formula \eqref{1.4} for the scalar case, the assumptions for $w_{1,\theta}$ in \eqref{4.1.2}, and formula \eqref{4.1.5} that 
\begin{equation}
    (\mathcal{S}(F_{\theta})(\gamma))_{11} = \prod_{y \in \gamma \cap \mathbb{Z}^d} w_{1,\theta}(y) = \mathcal{S}(w_{1,\theta})(\gamma).
\label{4.1.7}
\end{equation}
Next, using induction on $n$ - the number of factors in the product in \eqref{4.1.5}, we show that 
\begin{equation}
(\mathcal{S}(F_{\theta})(\gamma))_{12} = \mathcal{P}_Wu(\gamma),
\label{4.1.8}
\end{equation}
where $W$ is given by \eqref{4.1.4}: 
\\\\
\textbf{Case of n = 1:} Formula \eqref{4.1.8} holds trivially. \\\\
\textbf{Inductive Step:}
Let the statement that formulas \eqref{4.1.5}, \eqref{4.1.6} imply \eqref{4.1.8} be proved for a given $n \in \mathbb{N}.$ Then, for $n$ replaced by $n+1$ in \eqref{4.1.5}, we have that
\begin{equation}
    \mathcal{S}(F_\theta)(\gamma) = \begin{pmatrix}
        w_{1,\theta}(y_{n+1}) & w_{2,\theta}(y_{n+1})u(y_{n+1}) \\
        0 & 1
    \end{pmatrix} \times
    \atop
    \begin{pmatrix}
        \prod_{k = 1}^n w_{1,\theta}(y_k) & \sum_{k = 1}^{n}\bigg(w_{2,\theta}(y_k)\prod_{l = k+1}^nw_{1,\theta}(y_l)\bigg)u(y_k) \\
        0 & 1
    \end{pmatrix},
\label{4.1.9}
\end{equation}
and hence, we have 
\begin{equation}
\begin{aligned}
    \mathcal{S}(F_\theta)(\gamma)_{12} = \sum_{k = 1}^{n}\bigg(w_{2,\theta}(y_k)\prod_{l = k+1}^{n+1}w_{1,\theta}(y_l)\bigg)u(y_k) + w_{2,\theta}(y_{n+1})u(y_{n+1})
    \\
    = \sum_{k = 1}^{n+1}\bigg(w_{2,\theta}(y_k)\prod_{l = k+1}^{n+1}w_{1,\theta}(y_l)\bigg)u(y_k) = \mathcal{P}_Wu(\gamma),
\end{aligned}
\label{4.1.10}
\end{equation}
where $P_Wu(\gamma)$ is as in \eqref{4.1.6} with $n$ replaced by $n+1$. Hence, induction holds. This finishes the proof of Theorem 4.1.
\\\\
Let 
\begin{equation}
    S_{r,\theta} := \{z \in B_r \cap \mathbb{Z}^d: z + s\theta \notin B_r \cap \mathbb{Z}^d \ \forall s > 0 \}, \ r > 0, \ \theta \in \mathbb{S}^{d-1}.
\label{4.3}
\end{equation}
One can see that the structure of $S_{r,\theta}$ strongly depends on $\theta$: 
\begin{itemize}
    \item if $\theta \in \mathbb{S}^{d-1} \setminus \Omega$, then $S_{r,\theta} = B_r \cap \mathbb{Z}^d.$
    \item if $\theta$ is "very rational", for example, it coincides with one of the standard basis vectors in $\mathbb{R}^d$, then $S_{r,\theta}$ is a discrete analogue of the half sphere $\mathbb{S}^{d-1}_{r,\theta} := \{z \in B_r: z + s\theta \notin B_r  \ \forall s > 0 \}$.
\end{itemize} 
\textbf{Remark 3.2} 
\textit{Any $W = W(y,\theta)$ such that}
\begin{equation}
    W \textrm{ takes values in } \mathbb{C}\setminus\{0\}, \textrm{ supp}(W(\cdot,\theta) - 1) \subseteq B_r \cap \mathbb{Z}^d,
\label{4.2.1}
\end{equation}
\textit{admits presentation \eqref{4.1.4} with $w_{2,\theta} \equiv 1$ on $(B_r \cap\mathbb{Z}^d) \setminus S_{r,\theta}$.}
\\\\
Remark 3.2 follows from the following observations:
\begin{itemize}
    \item For $\gamma$ and $\{y_1,...,y_n\}$ in \eqref{4.1.5},\eqref{4.1.6}, we have that $\gamma \cap S_{r,\theta} = \{y_n\}$.
    \item As in \eqref{4.1.5}, we have that 
    \begin{equation}
        W(y_k,\hat{\gamma}) = w_{2,\theta}(y_k)\prod_{l = k+1}^{n}w_{1,\theta}(y_l), \ k = 1,...,n.
    \label{4.2.2}
    \end{equation}
    \item System \eqref{4.2.2} is uniquely solvable for $w_{2,\theta}(y_n), w_{1,\theta}(y_{n-1}),...,w_{1,\theta}(y_1)$ under the conditions that $W(y_k, \hat{\gamma}) \in \mathbb{C} \setminus \{0\}, \ w_{2,\theta}(y_j) = 1$ for $j = 1,...,n-1$.
\end{itemize}
Recall that in the continuous case as in \eqref{1.1},\eqref{1.2}, if
\begin{equation}
  f(x) = \begin{pmatrix}
      a(x) & u(x) \\
      0 & 0
  \end{pmatrix},
\label{4.4}
\end{equation}
then,
\begin{equation}
    \mathcal{S}(f)(\gamma) = \begin{pmatrix}
    \exp(-\mathcal{P}a(\gamma)) & -\mathcal{P}_Wu(\gamma) \\
    0 & 1
    \end{pmatrix},
\label{4.5}
\end{equation}
\begin{equation}
    W(y,\theta) = \exp\bigg(-\int_{0}^{\infty}a(y+s\theta)ds\bigg), \ y \in \mathbb{R}^d, \ \theta \in \mathbb{S}^{d-1},
\label{4.6}
\end{equation}
where $\mathcal{P}_W$ is the attenuated X-ray transform and $a$ is the attenuation coefficient; see \cite{N2002}, \cite{N2019}.\\
In this case, in the framework of discretization in formula \eqref{1.9}, we have that
\begin{equation}
    \exp(f(y_i)\Delta)  
    = \begin{pmatrix}
        w_{1,\Delta}(y_i) & w_{2,\Delta}(y_i)u(y_i) \\
        0 & 1
    \end{pmatrix},
\label{4.7}
\end{equation}
where 
\begin{equation}
    w_{1,\Delta}(y_i) = \exp(-a(y_i) \Delta), \ 
    w_{2,\Delta}(y_i) = \frac{\exp(-a(y_i)\Delta) - 1}{a(y_i)}.
\label{4.8}
\end{equation}
In addition, the step $\Delta$ between points of $\gamma \cap \mathbb{Z}^{d}, \ \gamma \in T'$, depends on the direction $\theta = \hat{\gamma}$. These observations explain the form of $F_\theta(z)$ in formula \eqref{4.1.1}.

\section{Continuous Case for Regularized Delta Functions}
We start with
\begin{equation}
    f(x) = \sum_{z \in \mathbb{Z}^d}w(z)f^{dis}(z)\chi_z(x-z), \ x \in \mathbb{R}^d,
\label{3.1}
\end{equation}
where 
\begin{equation}
f^{dis} \textrm{ takes values in } M(n,\mathbb{C}),  \ \textrm{supp}(f^{dis}
)\subseteq B_r \cap \mathbb{Z}^d, 
\label{3.2}
\end{equation}
\begin{equation}
    w > 0, \ \chi_z \textrm{ are characteristic functions of the balls }B_{\rho_z},
\label{3.3}
\end{equation}
where 
\begin{equation}
    \rho_z > 0,  \ 
    (\ B_{\rho_{z_1}} +\rho_{z_1})  \cap (B_{\rho_{z_2}} + \rho_{z_2}) = 0 \textrm{ if } z_1 \neq z_2, \ z,z_1,z_2 \in B_r \cap \mathbb{Z}^d.
\label{3.4}
\end{equation}
In fact, we consider $w(z)\chi_z(x-z)$ as approximations to $\delta(x-z)$ (or as regularized delta functions).\\
Let
\begin{equation}
T'' = \{\gamma \in T': z \in \gamma \cap (B_{\rho_z} + z) \textrm{ or } \gamma \cap (B_{\rho_z} + z) = \emptyset \ \forall z \in B_r \cap \mathbb{Z}^d\}.
\label{3.5}
\end{equation}
Note that 
\begin{equation}
    |T'_R \setminus T''| \rightarrow 0 \ \textrm{ as all } \rho_z \rightarrow 0, \ \forall R > 0,
\label{3.6}
\end{equation}
\begin{equation}
    T'_R := \bigcup_{\zeta \in B_R \cap \mathbb{Z}^d} T_\zeta \subset T',
\label{3.7}
\end{equation}
where $T_\zeta$ is defined as in \eqref{1.5}, $|\cdot|$ denotes the Lebesgue measure on $T_R'$, and $z \in B_r \cap \mathbb{Z}^d.$
\\\\
Let 
\begin{equation}
    F(z) = \exp(2\rho_zw(z)f^{dis}(z)), \ z \in B_r \cap \mathbb{Z}^d, \ F(z) = \mathbf{I}, \ z \in \mathbb{Z}^d \setminus \{B_r\}.
\label{3.8}
\end{equation}
\\\\
\textbf{Proposition 4.1: }\textit{Under assumptions \eqref{3.1}- \eqref{3.5}, and \eqref{3.8}, the following formula holds:}
\begin{equation}
    \mathcal{S}^{con}(f)(\gamma) = \mathcal{S}^{dis}(F)(\gamma), \  \ \gamma \in T'',
\label{3.1.1}
\end{equation}
where $\mathcal{S}^{con}$ is defined in \eqref{1.1},\eqref{1.2}, and $\mathcal{S}^{dis}$ is defined in \eqref{1.4}.
\\\\
\textbf{Proof of Proposition 4.1}
Using  \eqref{1.1}, \eqref{1.2}, under our assumptions, we get:
\begin{equation}
    \mathcal{S}^{con}(f)(\gamma) = \prod^{\longleftarrow}_{z \in \gamma \cap \mathbb{Z}^d \cap B_r}\exp(w(z)|\gamma \cap ( B_{\rho_z} + z)|f^{dis}(z)), \ \gamma \in T'',
    \label{3.1.2}
\end{equation}
where $|\cdot|$ denotes the length of the interval. 
Formula \eqref{3.1.1} follows from \eqref{3.1.2} and the observation that
\begin{equation}
|\gamma \cap (\ B_{\rho_z} + z)| = 2\rho_z, \  \gamma 
\in T'', \ z \in \gamma \cap \mathbb{Z}^d \cap B_r.
\end{equation}
This finishes the proof of Proposition 4.1. \\\\
In some cases, Proposition 4.1 implies that reconstructing $f$ from its continuous non-abelian X-ray transform $\mathcal{S}^{con}(f)$ reduces to reconstructing $F$ from its discrete non-abelian X-ray transform $\mathcal{S}^{dis}(F)$, under the assumption that $f$ is given by \eqref{3.1} with known $w$ and $\rho_z$. \\
A natural example in this connection consists of reconstructing $f$ from $\mathcal{S}(f)$ given on $\mathcal{T}^{\star}$ defined by \eqref{2.20} for the case when all $\rho_z$ in \eqref{3.3} are sufficiently small such that $\mathcal{T}^{\star} \subseteq T''$.

\section{Continuous Case for Piecewise Constant Functions}
Note that $\mathbb{R}^d$ can be written as:
\begin{equation}
    \mathbb{R}^d = \bigcup_{z \in \mathbb{Z}^d}\mathcal{U}_z,
\label{3.2.1}
\end{equation}
where
\begin{equation}
    \mathcal{U}_z = \mathcal{U} + z,\ \mathcal{U} = \{x = (x_1,\cdots, x_d) \in \mathbb{R}^d: -1/2 \leq x_j < 1/2, \ j = 1,\cdots,d\}.
\label{3.2.2}
\end{equation}
A regular function $f^{reg}$ on $\mathbb{R}^d$ supported in $B_r$ can be approximated by a piecewise constant function $f$ such that 
\begin{equation}
f(x) = 
\begin{cases}
    c(z), \ x \in \mathcal{U}_z, \  z  \in B_r \cap \mathbb{Z}^d,\\
    0, \ x \in \mathcal{U}_z, \  z  \in \mathbb{Z}^d \setminus B_r.
\end{cases}
\label{3.2.3}
\end{equation}
In particular, we have that 
\begin{equation}
    f(x) = \sum_{\zeta \in \mathbb{Z}^d} f^{dis}(\zeta)\chi_{\zeta}(x), \ x \in \mathbb{R}^d,
\label{3.2.4}
\end{equation}
where
\begin{equation}
    f^{dis}(z) = 
    \begin{cases}
        c(z), \ z \in B_r \cap \mathbb{Z}^d, \\
        0, \ z \in \mathbb{Z}^d \setminus B_r,
    \end{cases}
\label{3.2.5}
\end{equation}
and $\chi_\zeta$ is the characteristic function of $\mathcal{U}_{\zeta}.$ Here, $f^{reg}, f$, and $f^{dis}$ take values in $M(n, \mathbb{C}).$
\\\\
Let 
\begin{equation}
    \Sigma_\gamma = \{z \in \mathbb{Z}^d: \gamma \cap \mathcal{U}_z \neq \emptyset \}, \ \gamma \in T,
\label{3.2.7}
\end{equation}
where the points $\{z\}$ of $\Sigma_{\gamma}$ are ordered according to the order of intersection of $\mathcal{U}_z$ with $\gamma$.\\\\
If $f$ is defined as in \eqref{3.2.4}, then, in view of \eqref{1.1}, \eqref{1.2}, the continuous non-abelian X-ray transform is given by 
\begin{equation}
    \mathcal{S}^{con}(f)(\gamma) = \prod^{\longleftarrow}_{\zeta \in \Sigma_{\gamma}}\exp(|\gamma \cap \mathcal{U}_\zeta|f^{dis}(\zeta)), \ \gamma \in T,
\label{3.2.6}
\end{equation}
where $|\cdot|$ denotes the length of the interval. Motivated by formulas \eqref{1.4}, \eqref{3.2.6}, \eqref{3.2.n.1}, a discrete non-overdetermined version $\mathcal{S}^{\star}$ of continuous non-abelian X-ray  transform $\mathcal{S}$ in \eqref{1.1}, \eqref{1.2} can also be defined by the formula
\begin{equation}
    \mathcal{S}^{\star}(f^{dis})(\gamma) = \prod^{\longleftarrow}_{\zeta \in \mathbb{Z}^d}\exp(|\gamma \cap \mathcal{U}_\zeta|f^{dis}(\zeta)), \ \gamma \in \Gamma_r \subseteq T,
\label{3.2.8}
\end{equation}
where $f^{dis}$ is a $M(n,\mathbb{C})$-valued function on $\mathbb{Z}^d$ supported in $B_r \cap \mathbb{Z}^d$, and
\begin{equation}
    \#\Gamma_r = \#(B_r \cap \mathbb{Z}^d).
\label{3.2.n.1}
\end{equation}
When $f$ is defined as in \eqref{3.2.4}, it is obvious that 
\begin{equation}
    \mathcal{S}^{con}(f)(\gamma) = \mathcal{S}^{\star}(f^{dis})(\gamma), \textrm{ for } \gamma \in \Gamma_r.
\label{3.2.9}
\end{equation} 
We emphasize that $\mathcal{S}^{\star}$ defined in \eqref{3.2.8}, \eqref{3.2.n.1} is different from the discrete non-overdetermined X-ray transforms considered in Section 2. \\\\
In fact, the definition of a discrete version of continuous non-abelian X-ray transform via \eqref{3.2.8} is quite natural in the framework of numerical analysis. In connection with related definitions for the case of classical X-ray transform, see, for example, \cite{AS2004}. 
\\\\
The following important formula holds:
\begin{equation}
    \begin{aligned}
    \mathcal{S}^{\star}(f^{dis})(\gamma) = \bigg(\prod^{\longleftarrow}_{\zeta \in \Sigma_{\gamma,z}^+}\exp(|\gamma \cap \mathcal{U}_\zeta|f^{dis}(\zeta))\bigg)  \exp(|\gamma \cap \mathcal{U}_z|f^{dis}(z))
    \\ \times
\bigg(\prod^{\longleftarrow}_{\zeta \in \Sigma_{\gamma,z}^-}\exp(|\gamma \cap \mathcal{U}_\zeta|f^{dis}(\zeta))\bigg), \  z \in \gamma \in T, 
\end{aligned}
\label{3.2.10}
\end{equation}
where
\begin{equation}
    \Sigma_{\gamma,z}^{+} := \{\zeta \in \Sigma_\gamma: \gamma \textrm{ intersects } \mathcal{U}_\zeta \textrm{ after } \mathcal{U}_z \},
    \atop 
    \Sigma_{\gamma,z}^{-} := \{\zeta \in \Sigma_\gamma: \gamma \textrm{ intersects } \mathcal{U}_\zeta \textrm{ before } \mathcal{U}_z \}.
\label{3.2.n.2}
\end{equation}
Formula \eqref{3.2.10} follows from \eqref{3.2.8}.
\\
Using \eqref{3.2.10}, we suggest a layer-stripping reconstruction of $f^{dis}$ on $B_r \cap \mathbb{Z}^d$ from $\mathcal{S}^{\star}(f^{dis})$ on some appropriate $\Gamma_r$ in a similar way to the considerations in Subsection 2.2. This result is summarized in Theorem 5.5 below. One of the key points in this result consists of constructing $\Gamma_r$; see formulas \eqref{3.2.nn.8} for $d = 2$ and \eqref{3.2.nnn.1} for $d \geq 3$. 
\\\\
Let
\begin{equation}
    \mathring{\mathcal{U}}_z = \mathring{\mathcal{U}} + z,\ \mathring{\mathcal{U}} = \{x = (x_1, x_2) \in \mathbb{R}^2: -1/2 < x_j < 1/2, \ j = 1,2\}, \ z \in \mathbb{Z}^2.
\label{3.2.nn.84}
\end{equation}
Let $\gamma_y$ denote the oriented straight line such that 
    \begin{equation}
        y \in \gamma_y, \ \hat{\gamma}_y = \frac{(-y_2,y_1)}{|y|}, \ y \in \mathbb{R}^2 \setminus \{0\}, \atop\hat{\gamma}_y = \theta,  \ y = 0, \textrm{ for some fixed } \theta \in \mathbb{S}^1.
    \label{3.2.nn.5}
    \end{equation}
We assume that 
\begin{equation}
    \sup_{z \in B_r \cap \mathbb{Z}^d}\|f^{dis}(z)\| \leq M \textrm{ for some } M > 0,
\label{3.2.nnn.4}
\end{equation}
where
\begin{equation}
\|f^{dis}(z)\| := \sqrt{\textrm{Trace}(f^{dis}(z)^{\dagger}f^{dis}(z))},
\label{3.2.nnn.5}
\end{equation} 
where $f^{dis}(z)^{\dagger}$ denotes the conjugate transpose. \\
We recurrently define the domains $D_i$, and the sets $V_i$, $Z_i,$ and $W_i,  \ i = 1,...,J$:
\begin{itemize}
    \item For $i =1$:
    \begin{equation}
    D_1 := \bigcup_{z \in B_r \cap \mathbb{Z}^2} \overline{{\mathcal{U}}}_z,
    \label{3.2.n.84}
    \end{equation}
    \begin{equation}
        V_1 := \{x \in D_1 : |x| \geq |y|, \ \forall y \in D_1 \},
    \label{3.2.n.3}
    \end{equation}
    \begin{equation}
        Z_1 := \{z \in \mathbb{Z}^2: \exists x = x(z) \in V_1 \textrm{ such that } x \in \overline{\mathcal{U}}_z \}.
    \label{3.2.n.4}
    \end{equation}
    Note that choice of $x = x(z)$ in \eqref{3.2.n.4} may not be unique, but is assumed to be fixed. In any case, each $x$ in $V_1$ belongs to just one $\overline{\mathcal{U}}_z$ in \eqref{3.2.n.84} and is one of its vertices.  Hence, $z = z(x)$ is well defined for $x \in V_1$. See also Remark 5.1 below.\\
    In addition:
    \begin{equation}
        W_1  := \{y \in \mathbb{R}^2: y = y(z) \in \mathring{\mathcal{U}}_z, \ z \in Z_1, \ \gamma_{y(z)} \cap \ \overline{(D_1 \setminus \overline{\mathcal{U}}_z)} = \emptyset \},
    \label{3.2.nn.6}
    \end{equation}
    where $\mathring{\mathcal{U}}_z, \ \gamma_y$ are defined in \eqref{3.2.nn.84}, \eqref{3.2.nn.5}, respectively. Note that there are uncountably many choices of $y(z)$ in \eqref{3.2.nn.6}, but the choice is assumed to be fixed. In addition, it is convenient to assume that 
    \begin{equation}
        y(z) = \lambda(z)x(z) \textrm{ for some } \lambda(z) \in (0,1),
    \label{3.2.nn.9}
    \end{equation}
    where $x(z)$ is as in \eqref{3.2.n.4}; in this case, $\hat{\gamma}_{x(z)} = \hat{\gamma}_{y(z)}$. See also Remarks 5.2, 5.3 below.
    \item The inductive step is as follows:
    \begin{equation}
        D_i := \overline{D_{i-1} \setminus\bigcup_{z \in Z_{i-1}} \overline{{\mathcal{U}}}_z} = \bigcup_{z \in (B_r \cap \mathbb{Z}^2) \setminus \cup_{j = 1}^{i-1}Z_j}\overline{\mathcal{U}}_z, \ i = 2,...,J,
    \label{3.2.n.6}
    \end{equation}
    \begin{equation}
       V_i := \{x \in D_i : |x| \geq |y|, \ \forall y \in D_i \}, \ i = 2,...,J,
    \label{3.2.n.7}
    \end{equation}
    \begin{equation}
        Z_i := \{z \in \mathbb{Z}^2: \exists x = x(z) \in V_i \textrm{ such that } x \in \overline{\mathcal{U}}_z \}, \ i = 2,...,J,
    \label{3.2.n.8}
    \end{equation}
    \begin{equation}
    \begin{aligned}
         W_i  := \{y \in \mathbb{R}^2: y = y(z) \in \mathring{\mathcal{U}}_z, \  z \in Z_i, \ \gamma_{y(z)} \cap \ \overline{(D_i \setminus \overline{\mathcal{U}}_z)} = \emptyset \}, \atop i = 2,...,J,
    \label{3.2.nn.7}
    \end{aligned}
    \end{equation}
    where 
    \begin{equation}
     D_J := \overline{\mathcal{U}}_0, \ Z_J = \{0\}.
     \label{3.2.nn.1}
    \end{equation}
    As in formulas \eqref{3.2.n.84} - \eqref{3.2.n.4}, choice of $x = x(z)$ in \eqref{3.2.n.8} may not be unique, but is assumed to be fixed. In any case, each $x$ in $V_i$ belongs to just one $\overline{\mathcal{U}}_z$ in the right hand side of \eqref{3.2.n.6} and is one of its vertices. Hence, $z = z(x)$ is well defined for $x \in V_i$. See also Remark 5.1 below.\\ 
    As in \eqref{3.2.nn.6}, there are uncountably many choices of $y(z)$ in \eqref{3.2.nn.7}, but the choice is assumed to be fixed. In addition, it is convenient to assume that \eqref{3.2.nn.9} holds, where $x(z)$ is as in \eqref{3.2.n.8}. See also Remarks 5.2, 5.3 below.
    \end{itemize}
Proceeding from the sets $W_i$, we define
\begin{equation}
    \Gamma_{r,i} := \bigcup_{y \in W_i}\gamma_y, \ i \in \{1,...,J\},
\label{3.2.nnn.3}
\end{equation}
\begin{equation}
    \Gamma_r := \bigcup_{y \in W}\gamma_y = \bigcup_{z \in B_r \cap \mathbb{Z}^2} \gamma_{y(z)},
\label{3.2.nn.8}
\end{equation}
where
\begin{equation}
    \\W = \bigcup_{i = 1}^J W_i,
\label{3.2.nn.4}
\end{equation}
and $y(z)$ is as in \eqref{3.2.nn.6} and \eqref{3.2.nn.7}. \\ By construction of $W_i$ and $Z_i$, we also have that
\begin{equation}
   \# W = \#(\bigcup_{i = 1}^JZ_i) = \#(B_r \cap \mathbb{Z}^2) =: N_r,
\label{3.2.nn.2}
\end{equation} 
\begin{equation}
    \bigcup_{i = 1}^{J}Z_i = B_r \cap \mathbb{Z}^2.
\label{3.2.nn.3}
\end{equation}
\\\\
\textbf{Remark 5.1} 
\textit{The property that each $x \in V_i$ is a vertice of one of the squares $\overline{\mathcal{U}}_z$ arising in \eqref{3.2.n.84} and in the right hand side of \eqref{3.2.n.6} follows from the observation that for each closed interval of a line (edge of a closed square in our case), the maximal distance from the origin is attained at one of its boundary points. A point of $x \in V_i$ cannot be a shared vertice of two squares $\overline{\mathcal{U}}_{z_1}$ and $\overline{\mathcal{U}}_{z_2}$, since the line segment $\overline{y_1y_2}$ formed by their edges connected at $x$ must attain its maximal distance from the origin at either $y_1$ or $y_2$. In addition, non-uniqueness of $x(z)$ in \eqref{3.2.n.4}, \eqref{3.2.n.8} is illustrated by $\overline{\mathcal{U}}_0$, where each of its four vertices are at the same distance from the origin.}
\\\\
\textbf{Remark 5.2} \textit{In connection with the definition of $W_i$ in \eqref{3.2.nn.6}, \eqref{3.2.nn.7}, it is important to note that }
\begin{equation}
    \tilde{\mathcal{U}}_z \neq \emptyset, \ z \in Z_i, \ i = 1,...,J,
    \label{3.2.nn.10}
\end{equation}
\textit{where}
\begin{equation}
    \tilde{\mathcal{U}}_z := \{y \in \mathring{\mathcal{U}}_z:  \gamma_{y} \cap \ \overline{(D_i \setminus \overline{\mathcal{U}}_z)} = \emptyset\}, \ z \in Z_i.
    \label{3.2.nn.11}
\end{equation}
\textit{The reason for \eqref{3.2.nn.10} is that}\\
\begin{equation}
    \gamma_x \cap D_i = \{x\}, \ \textrm{dist}(\gamma_x, \overline{(D_i \setminus \overline{\mathcal{U}}_z)}) > 0, \ \forall x \in V_i, \ i = 1,...,J.
\label{3.2.nnn.9}
\end{equation}
\textit{This is based on the fact that $x$, by definition of $V_i$, is an extremal point of $D_i$, the orthogonality $\langle y - x, x \rangle = 0, \ \forall y \in \gamma_x$, and the compactness of $\overline{(D_i \setminus \overline{\mathcal{U}}_z)}$. In particular, if $y \in \gamma_x \cap D_i, \ y \neq x$, then $|y| > |x|,$ contradicting the extremality of $x$ in $D_i$. One can see that $W_i \subset \bigcup_{z \in Z_i} \tilde{\mathcal{U}}_z.$}
\\\\
\textbf{Remark 5.3}
\textit{In \eqref{3.2.nn.9}, we use that for $x(z) = (x_1(z), x_2(z))$, $|x_k(z)| > |y_k|, \ k  = 1, 2,$ for all $y = (y_1,y_2) \in \mathring{\mathcal{U}}_z$. Hence, in particular,} 
\begin{equation}
    \exists \varepsilon_0(z)  \in (0,1) \textrm{ such that } (1 - \varepsilon)x(z) \in \mathring{\mathcal{U}}_z, \ \forall \varepsilon \in (0, \varepsilon_0(z)), \ z \in B_r \cap \mathbb{Z}^2.
\label{3.2.nnn.10}
\end{equation}
\\
\textbf{Remark 5.4}
\textit{In view of \eqref{3.2.nnn.9}, \eqref{3.2.nnn.10}, we have that}
\begin{equation}
    0< |\gamma_{y(z)} \cap \mathcal{U}_z| < (\log 2)/M, \ \textrm{dist}(\gamma_{y(z)}, \overline{(D_i \setminus \overline{\mathcal{U}}_z)}) > 0,\  z \in Z_i,\ i = 1,...J, 
\label{3.2.nnn.11}
\end{equation}
\textit{if $y(z) = (1 - \varepsilon)x(z)$ for sufficiently small $\varepsilon > 0$, and for arbitrary fixed $M > 0$.}
\\\\
The simplest construction of $\Gamma_r$ for $d \geq 3$ consists of repeating the two dimensional construction on each two-dimensional plane $\Xi = \Xi_{z_3,...,z_d}$ defined in \eqref{2.n.1}. More precisely, 
\begin{equation}
    \Gamma_r : = \bigcup_{z_3^2 +...+z_d^2 \leq r^2} \Gamma_{r, (z_3,...,z_d)}, \ z_3, ...,z_d \in \mathbb{Z},
\label{3.2.nnn.1}
\end{equation}
where $\Gamma_{r,(z_3,...z_d)}$ is constructed on $\Xi_{z_3,...z_d}$ as in \eqref{3.2.nn.8}, \eqref{3.2.nnn.11} with $r$ replaced by $\sqrt{r^2 - z_3^2 - ... - z_d^2}.$
\\\\
\textbf{Theorem 5.5} \textit{Let $f^{dis}$ be a $M(n, \mathbb{C})$-valued function on $\mathbb{Z}^d$, satisfying \eqref{3.2.5}, \eqref{3.2.nnn.4}. Then: 
\begin{itemize}
    \item $\mathcal{S}^{\star}(f^{dis})$ on $\Gamma_r$ uniquely determines $f^{dis}$, where $\Gamma_r$ is defined
by \eqref{3.2.nn.8}, \eqref{3.2.nnn.11} for $d = 2$, and by \eqref{3.2.nnn.1} for $d \geq 3$.
\item This determination is given via the two-dimensional formulas \eqref{3.2.n.15}, \eqref{3.2.n.16} for $d = 2$, which are adopted for each slice $\Xi$ with fixed $z_3,...,z_d$ for $d \geq 3$.
\item For $d = 2$, formulas \eqref{3.2.n.15}, \eqref{3.2.n.16} also yield a reconstruction of $f^{dis}$ on $\bigcup_{i=1}^jZ_i$ from $\mathcal{S}^{\star}(f^{dis})$ on $\bigcup_{i=1}^j\Gamma_{r,i}, \ j = 1,2,...J$, where $Z_i$ and $\Gamma_{r,i}$ are defined by \eqref{3.2.n.4},\eqref{3.2.n.8}, and \eqref{3.2.nnn.3}, respectively.
\end{itemize}}
\textbf{Proof of Theorem 5.5} 
We first consider the case of $d = 2$. \\
By assumption \eqref{3.2.nnn.4} and \eqref{3.2.nnn.11}, we have 
\begin{equation}
    \|\Lambda(z,z) f^{dis}(z)\| < \log 2, \ \forall z \in B_r \cap \mathbb{Z}^2, 
\label{3.2.nnn.6}
\end{equation}
where
\begin{equation}
    \Lambda(z,\zeta) := |\gamma_{y(z)} \cap \mathcal{U}_\zeta|, \ z,\zeta \in B_r \cap \mathbb{Z}^2,
\label{3.2.nnn.12}
\end{equation}
and the norm $\|\cdot\|$ is defined in \eqref{3.2.nnn.5}. Hence, 
\begin{equation}
    \log(\exp(\Lambda(z,z)f^{dis}(z)) = \Lambda (z,z)f^{dis}(z).
\label{3.2.nnn.8}
\end{equation}
where 
\begin{equation}
    \log(\mathbf{I} + \Delta) = \sum_{m = 1}^{\infty}(-1)^{m+1}\frac{\Delta^m}{m}, \ \|\Delta\| < 1.
\label{3.2.nnn.7}
\end{equation}
The convergence of the series in \eqref{3.2.nnn.7} and the identity \eqref{3.2.nnn.8} follow from the condition \eqref{3.2.nnn.6} and standard results on matrix logarithm; see, for example, Theorem 2.8 in \cite{H2015}.
\\
Next, using formulas \eqref{3.2.10}, \eqref{3.2.n.4},
\eqref{3.2.nn.6}, \eqref{3.2.n.8}, \eqref{3.2.nn.7}, \eqref{3.2.nnn.11}, we get 
\begin{equation}
    \exp(|\gamma_{y(z)} \cap \mathcal{U}_z|f^{dis}(z)) = \mathcal{S}^{\star}(f^{dis})(\gamma_{y(z)}),  \ z \in Z_1,
\label{3.2.n.15}
\end{equation}
\begin{equation}
    \begin{aligned}
    \exp(|\gamma_{y(z)} \cap \mathcal{U}_z|f^{dis}(z)) = \bigg(\prod^{\longleftarrow}_{\zeta \in (\cup_{k = 1}^{i}Z_k) \cap \Sigma_{\gamma_{y(z)},z}^+}\exp(|\gamma_{y(z)} \cap \mathcal{U}_\zeta|f^{dis}(\zeta))\bigg)^{-1} \mathcal{S}^{\star}(f^{dis})(\gamma_{y(z)})
    \\\times
\bigg(\prod^{\longleftarrow}_{\zeta \in (\cup_{k = 1}^{i}Z_k) \cap \Sigma_{\gamma_{y(z)},z}^-}\exp(|\gamma_{y(z)} \cap \mathcal{U}_\zeta|f^{dis}(\zeta))\bigg)^{-1}, \ z \in Z_{i+1},  \ i = 1,...,j-1.
\label{3.2.n.16}
\end{aligned}
\end{equation} 
In turn, using \eqref{3.2.nnn.8}, \eqref{3.2.nnn.7}, \eqref{3.2.n.15}, \eqref{3.2.n.16}, we get
\begin{equation}
    f^{dis}(z) = \frac{1}{|\gamma_{y(z)} \cap \mathcal{U}_z|}\log\{\mathcal{S}^{\star}(f^{dis})(\gamma_{y(z)})\},  \ z \in Z_1,
\label{3.2.nnn.13}
\end{equation}
\begin{equation}
\begin{aligned}
    f^{dis}(z) = \frac{1}{|\gamma_{y(z)} \cap \mathcal{U}_z|} \log \bigg\{\bigg(\prod^{\longleftarrow}_{\zeta \in (\cup_{k = 1}^{i}Z_k) \cap \Sigma_{\gamma_{y(z)},z}^+}\exp(|\gamma_{y(z)} \cap \mathcal{U}_\zeta|f^{dis}(\zeta))\bigg)^{-1} \mathcal{S}^{\star}(f^{dis})(\gamma_{y(z)})
    \\\times
\bigg(\prod^{\longleftarrow}_{\zeta \in (\cup_{k = 1}^{i}Z_k) \cap \Sigma_{\gamma_{y(z)},z}^-}\exp(|\gamma_{y(z)} \cap \mathcal{U}_\zeta|f^{dis}(\zeta))\bigg)^{-1}\bigg\}, \ z \in Z_{i+1},  \ i = 1,...,j-1,
\end{aligned}
\label{3.2.nnn.14}
\end{equation}
Recall that $\Sigma^+_{\gamma,z},\Sigma^-_{\gamma,z}$ are defined by \eqref{3.2.n.2}. \\
Formulas \eqref{3.2.nnn.13}, \eqref{3.2.nnn.14} yield a layer-stripping reconstruction of $f^{dis}$ on $\bigcup_{i = 1}^{j}Z_i$ from $\mathcal{S}^{\star}(f^{dis})$ on $\bigcup_{i = 1}^{j}\Gamma_{r,i}, \ j = 1,...,J$. \\
In particular, when $j = J$, \eqref{3.2.nnn.13}, \eqref{3.2.nnn.14} give the reconstruction formulas of $f^{dis}$ on $B_r \cap \mathbb{Z}^d$ from $\mathcal{S}^{\star}(f^{dis})$ on $\Gamma_r$. This finishes the proof for $d = 2$. \\
In dimension $d \geq 3$, we adopt the reconstruction formulas \eqref{3.2.nnn.13}, \eqref{3.2.nnn.14} for each two-dimensional plane $\Xi_{z_3,...,z_d}$ defined by \eqref{2.n.1}, where $z_3,...,z_d \in \mathbb{Z}, \ z_3^2 +...+z_d^2 \leq r^2,$ in a similar way to the proof of Theorem 2.1. This completes the proof of Theorem 5.5.
\\\\
\textbf{Remark 5.6}\textit{ Assumption \eqref{3.2.nnn.11} for $\Gamma_r$ is very essential for Theorem 5.5, in general. The simplest counterexample to Theorem 5.5 without assumption \eqref{3.2.nnn.11} is as follows: }\\
\textit{Let $d = 2, \ n = 1, \ 0 < r < 1.$ Then,}
\begin{equation}
    \mathcal{S}^{\star}(f^{dis})(\gamma_{y(0)}) = \exp(\Lambda(0,0)) \textrm{ for } f^{dis}(0) = 1 + \frac{(2k\pi i)}{\Lambda(0,0)}, \  \forall k \in \mathbb{Z}, 
\label{3.2.nnn.15}
\end{equation}
\textit{where $\Lambda(0,0)$ is defined by \eqref{3.2.nnn.12} and $\Gamma_r = \{\gamma_{y(0)}\}$ ($\#\Gamma_r = 1$). In particular, one can see that the example of  non-uniqueness in \eqref{3.2.nnn.15} is given for reconstructing $f^{dis}$ in \eqref{3.2.4}, \eqref{3.2.5} from its $\mathcal{S}^{\star}(f^{dis})$ on $\Gamma_r$ already for the scalar case ($n = 1$) related to the classical X-ray tomography with $f^{dis}$ taking values in $\mathbb{C}.$ }
\\\\
\textbf{Remark 5.7}\textit{ Assumption \eqref{3.2.nnn.11} for $\Gamma_r$ in Theorem 5.5 can be replaced by additional a priori assumptions on $f^{dis}$. For example, for the case when $n = 1$ and $f^{dis}$ is real-valued (basic case of the classical X-ray tomography), Theorem 5.5 holds without assumption \eqref{3.2.nnn.11}.}
\\\\
\textbf{Remark 5.8}
\textit{For the case when $n = 1$ and $f^{dis}$ is real-valued, Theorem 5.5 continues studies in \cite{NS2025} for reconstructing $f^{dis}$ on $B_r \cap \mathbb{Z}^d$ from data equivalent to $\mathcal{S}^{\star}(f^{dis})$ on an appropriate $\Gamma_r$ such that $\#\Gamma_r = \#(B_r \cap \mathbb{Z}^d)$. In this connection, \cite{NS2025} suggests an approximate reconstruction for the case when $\Gamma_r = \mathcal{T}^{\star}_{0,r}$, where $\mathcal{T}^{\star}_{\alpha, \beta}$ is defined in \eqref{2.25}. In particular, in the present work, we succeeded to find an exact version of this approximate reconstruction by making small modifications to $\Gamma_r$.}
\\\\
\textbf{Remark 5.9}
\textit{For practical applications, $\mathcal{S}^{con}(f) = \mathcal{S}^{\star}(f^{dis})$ (as in formula \eqref{3.2.9}) on $\Gamma_r$ considered in Theorem 5.5 can be found via interpolations from $\mathcal{S}^{con}(f)$ given on more standard family of oriented straight lines.}
\\\\
\textbf{Remark 5.10}
\textit{In dimension $d \geq 3$, layer-stripping reconstructions for the continuous non-abelian X-ray transforms go back to \cite{N2002}. In dimension $d = 2$, to our knowledge, results on layer-stripping reconstruction for the continuous non-abelian X-ray transform were not given in the literature before Theorem 5.5 of the present work. For the case of the classical continuous X-ray transform, such results go back to \cite{C1963}. In connection with uniqueness and non-uniqueness results for the case of layer-stripping reconstruction for weighted continuous X-ray transform, see \cite{F1986}, \cite{BQ1987}, \cite{B1993}, \cite{GN2019}.}

\section*{Acknowledgement}
This work was fulfilled during the internship of the first author in the Centre de Mathématiques Appliquées of École polytechnique in June-August 2025 in the framework of research program for international talents.

\bigskip
Pranav Gupta \\
Department of Mathematics, \\
National University of Singapore, \\
21 Lower Kent Ridge Rd, 119077, Singapore \\
Email: pranav.gupta@u.nus.edu
\\\\
Roman G. Novikov\\
CNRS, CMAP, École polytechnique, \\
Institut Polytechnique de Paris, 91120 Palaiseau, France\\
E-mail: roman.novikov@polytechnique.edu

\end{document}